\documentstyle[12pt]{article}
\newcommand{\braid}[2]{{#1}$\lower4pt\hbox{$\oo\atop\raise4pt
           \hbox{$\scriptscriptstyle\Phi $}$}${#2}}
\newcommand{\twist}[2]{{#1}${\,\scriptscriptstyle \Phi}\atop\raise9pt
           \hbox{$\scriptstyle\oo$} ${#2}}
\newcommand{\BX}{\lower2pt\hbox{$\Box$}}

\newcommand{\D}{\partial}

\newcommand{\tz}{\tilde{z}}
\newcommand{\tp}{\tilde{p}}
\newcommand{\tq}{\tilde{q}}
\newcommand{\tu}{\tilde{u}}

\newcommand{\tv}{\tilde{v}}
\newcommand{\td}{\tilde{\partial}}
\newcommand{\tgm}{\tilde{\gamma}}

\newcommand{\ve}{\varepsilon}

\newcommand{\be}{\begin{eqnarray}}
\newcommand{\ee}{\end{eqnarray}}
\newcommand{\n}{\nonumber }
\newcommand{\oo}{\otimes}
\newcommand{\al}{\alpha}
\newcommand{\bt}{\beta}
\newcommand{\la}{\lambda}
\newcommand{\si}{\sigma}
\newcommand{\gm}{\gamma}

\begin{document}
\begin{titlepage}
\hspace{4.7in}{LPTHE 98-57}
\vspace{0.3in}
\begin{center}
{\Large \bf Twist-related geometries on q-Minkowski space}
\end{center}
\vspace{0.3cm}
\begin{center}
 {\Large P. P. Kulish\footnote{Partially supported by the RFFI grant 98-01-00310.}}
\end{center}
\vspace{0.2cm}
\begin{center}
St.Petersburg Department of the Steklov
Mathematical Institute,\\
Fontanka 27, St.Petersburg, 191011, Russia; \\
Laboratoire de Physique Th\'eorique et Hautes \'Energies,\\
Universit\'e Pierre et Marie Curie (Paris VI), 75252, France\\
(kulish@pdmi.ras.ru)
\end{center}
\vspace{0.2cm}
\begin{center}
{\Large A.I. Mudrov}
\end{center}
\begin{center}
\vspace{0.2in}
 Department of Theoretical Physics,
Institute of Physics, St.Petersburg State University,
St.Petergof, St.Petersburg, 198904, Russia\\
(aimudrov@DG2062.spb.edu)
\end{center}
\vspace {0.3in}
\begin{center}
  { Abstract}
\end{center}
{
The role of the quantum universal enveloping algebras of symmetries in
constructing non-commutative geometry of the space-time
including vector bundles, measure, equations of motion and their
solutions is discussed. In the framework of the twist theory the
Klein-Gordon-Fock and Dirac equations on the
quantum Minkowski space are studied from this point of view for the
simplest quantum  deformation of the Lorentz algebra induced by its
Cartan subalgebra twist.
}
\end{titlepage}

\section{Introduction}

The problem of building non-commutative models of the space-time
arose in connection with attempts to formulate quantum field
theory on this new basis, that could be free from difficulties
inherent to the classical differential geometry description.
Numerous publications on that matter address it in various aspects, using
deformation of the group of relativistic symmetries as the starting
point for further considerations. The first step on that way
is to define the  quantum group and its comodules along the line of FRT
method \cite{RTF}. At the next stage differential calculus consistent
with the group coaction is constructed, following
\cite{W3,ShirWessZum,Schir}, and the  equations of motion and
observable algebra representations are studied
\cite{CD,LRTN,AzKulRod,PilSchmWes,Kul,AzKulRod1}. The programme
just outlined thus gives  the priority to the concept of group in
proceeding to non-commutative geometry of the space time.
In the present paper such a transition is performed on the basis
of  twist \cite{D1,D4} of a universal enveloping Lie algebra
as the algebra of symmetries of a geometrical space.
Twist preserves algebraic structure of a quantum algebra but changes
its comultiplication. Twist is accompanied by deformation of the
product in the  algebra of functions on quantum group, and that process
is entirely formulated in terms of the universal twisting 2-cocycle.
This observation can be put into the foundation of a systematic
approach to study whole classes of geometries corresponding to
twist-related Hopf algebras. Quantum deformations of the Lorentz algebra are
known to possess rich twist structure \cite{M3}, and there are explicitly
built the twisting cocycles for certain quantum-deformed Poincar{\'e} algebras
\cite{LRTN,M7}, so it is natural to apply this machinery to
investigate properties of the deformed space-time.

Recall some principal facts concerning the subject under consideration.
Let ${\cal H}$ be a Hopf algebra possessing an invertible element (twist)
$$\Phi =\sum_i \Phi_i^{(1)} \otimes \Phi_i^{(2)}= \Phi_1 \otimes \Phi_2
\in{\cal H}\otimes {\cal H},$$
which satisfies the equation  \cite{D4}
\be
  (\Delta\oo id)(\Phi) \Phi_{12} & = & (id\oo\Delta)(\Phi) \Phi_{23}
  \label{TE}
\ee
and the normalizing condition
\be
(\ve\oo id)(\Phi) = 1 = (id \oo\ve)(\Phi).
 \label{NC}
\ee
It allows to define a Hopf algebra  $\tilde {\cal H}$,
which coincides with ${\cal H}$ as an associative algebra but has
the coproduct
\be
\tilde\Delta(h) &=& \Phi^{-1} \Delta(h) \Phi.\label{TCP}
\ee
The counit remains unchanged due to (\ref{NC}).
The new antipode differs from the old one via the similarity
transformation $\tilde S(h)\equiv  u S(h) u^{-1}$ where
$ u=(\Phi^{-1}_{1}) S(\Phi^{-1}_{2})$.
Provided algebra ${\cal H}$ is quasitriangular, so will be
$\tilde {\cal H}$ and the new universal  R-matrix will be expressed
through the old one $R$ and the twisting cocycle:
$\tilde R = \Phi^{-1}_{21} R \Phi$.
Twist establishes an equivalence relation among Hopf algebras which
is not an isomorphism in the common sense as the twisted
coproduct may differ significantly from original, e. g. be
no more cocommutative. Twist-equivalence, having all the necessary
features (symmetry, reflexivity, and transitivity), manifests itself
in the equivalence of the monoidal categories of Hopf algebra
representations, and also in relations between objects of the
corresponding geometrical spaces, that is the subject of
the present paper.

A specific case of the twist construction very much involved in the theory
of quantum Lorentz algebra deformations is a twisted
tensor product of two Hopf algebras \cite{RSTS}. Let ${\cal H}$
be the ordinary tensor product ${\cal A}\oo {\cal B}$, and an element
$\Phi\in {\cal A}\oo {\cal B}$
$\sim {\cal A}\oo 1 \oo 1\oo {\cal B}
\subset{\cal H}\otimes {\cal H}$  satisfies the factorizability conditions
\be
 \Delta_{\cal A}(\Phi) =
 \Phi_{13}\Phi_{23}\in {\cal A}\oo{\cal A}\oo{\cal B}, & \quad &
 \Delta_{\cal B}(\Phi)=
 \Phi_{13}\Phi_{12}\in {\cal A}\oo{\cal B}\oo{\cal B}   \label{TTP}
\ee
Then $\Phi$ is a solution to the equation (\ref{TE}), and in this case
twist $\tilde {\cal H}$ of the original algebra with $\Phi$ will be
denoted \twist{${\cal A}$}{${\cal B}$}.

Recall that the $*$-operation (real form) on a Hopf algebra
${\cal H}$ means antilinear involutive algebra anti-automorphism
and coalgebra automorphism. Due to the identity $S*=*S^{-1}$,
which always holds due to the uniqueness of the antipode,
one can re-define the real form as $S^{2n}*$ for any integer number $n$.
We can also consider homomorphic and anti-cohomomorphic antilinear
operations of the kind $\theta=S^{2n+1}*$. Necessity of introducing
them is accounted for the following.
If some ${\cal H}$-comodule algebra ${\cal A}$ with
coaction  $\bt\colon {\cal A}\to {\cal H}\oo {\cal A} $  is endowed
with anti-involution $a\to \bar a$ consistent with $\bt$, it
 must satisfy the equality
 $\bt( \bar a) = (*\oo \bar{ }\> )\bt(a)$.
 In this case ${\cal H}$ plays the role of the function algebra on
a quantum group of transformation of some manifold. The
role of the universal enveloping algebra ${\cal H}^*$  is different
from the geometrical point of view because it
does act on ${\cal A}$. To ensure consistency between
the real forms and the action one has to require
$\overline{(h a)} = S(h^*) \bar a$, $h\in {\cal H}^*, a \in {\cal A}$.
So, further on we will mean by the real form of a quantum
algebra a homomorphic and anti-cohomomorphic antilinear
involution $\theta = S\circ *$.

It is important for the further geometrical considerations
to investigate how real forms survive twist.
From the formula
(\ref{TCP}) it follows that, condition
\be
(\theta\oo \theta )(\Phi)=\tau(\Phi),
 \label{CTE0}
\ee
fulfilled, the same involution is defined on $\tilde{\cal H}$.
For the $*$-operation the analogous natural requirement is
\be
 \Phi^*=\Phi^{-1}.
 \label{CTE}
\ee
However, condition  (\ref{CTE0}) is not always the case.
Take for example the jordanian quantization of the Borel
$sl(2)$ subalgebra. It is generated by two elements
$H$ and $X^-$ subject to relation  $[H,X^-]=-2 X^-$,
and the coproduct acquires a more simple form after
transition to the generator $\si=-\frac{1}{2}\ln(1-2X^-)$:
$\tilde\Delta(\si)=\si\oo 1 + 1\oo \si$,
$\tilde \Delta(H)=H\oo e^{2\si}+1\oo H$.
The deformation is induced by twist with the element
$\Phi=e^{- H\oo \si}$. It is easy to see that
the classical $*$-operation given by correspondence
$H^*=-H$, $(X^-)^*=X^-$ can be brought through the twist because of
(\ref{CTE}).
But for $\theta=S*$  one has
$(\theta\oo \theta)(\Phi)=\Phi^{-1}\not= \tau(\Phi)$,
thus the same $\theta$ cannot be defined on the twisted algebra.
Yet such a real form does exist and one can see that specifying it as
$\si\to -\si$, $H\to H e^{-2\si}$. The problem consists in
establishing a correspondence between this involution and
classical $\theta$. It is clear that we must find an automorphism
transforming the twisting element in such a way that it should
contain the factor $\tau(\Phi)$. The problem is typical for the class
of jordanian-like quantum algebras studied in \cite{M7,M,KLM}.
The key role in resolving it belongs to the identity
$(u\oo u)\tau (S\oo S)(\Phi)=\Phi^{-1}\Delta(u)$ fulfilled
for any solution to the twist equation \cite{D4}
(the element $u$ is exactly the same like that
taking part in the definition of the twisted antipode).
Let element $\Phi$ and  $*$-operation obey (\ref{CTE}).
Then we have
\be
\tau(\theta\oo\theta)(\Phi)=\tau(S\oo S)(\Phi^*)=
\tau(S\oo S)(\Phi^{-1})=\Delta(u^{-1})(\Phi)(u\oo u).
\label{AuxE}
\ee
The direct calculation shows that the mapping
$\tilde\theta$ given by the formula $\tilde\theta(h)=u\theta(h)u^{-1}$,
is anti-cohomomorphic:
$$(\tilde\theta\oo\tilde\theta)(\tilde\Delta(h))= \tau\tilde\Delta(\tilde\theta(h)).$$
Thus $\tilde\theta$ is an antilinear involutive algebraic  and
anti-coalgebraic automorphism of the algebra $\tilde{\cal H}$.

Quantizations of the Lorentz algebra  $U(so(1,3))$ are in close relations
with quantum algebra $U_q(sl(2))$,
so we remind here the basic facts of that theory.

The standard Drinfeld-Jimbo solution $U_q(sl(2))$ is build on the
generators $H$ and $X^{\pm}$ such that
$$[H, X^{\pm}]  =  \pm 2 X^{\pm}, \quad  [X^+ , X^-] =
\frac{{\rm sh}(\la H)}{{\rm sh}(\la)},$$
$$\Delta (H) = 1\oo H + H\oo 1, \quad
\Delta (X^{\pm}) = e^{-\frac{\la}{2} H}\oo X^{\pm} +
                 X^{\pm}\oo  e^{\frac{\la}{2} H},$$
where we set  $q=e^\la$.
The universal R-matrix is \cite{D0}
\be
  {\cal R}_q & = &
  e^{\frac{\la}{2} H \oo H}
  \sum_{n=0}^{\infty}\frac{(1-q^{-2})^n}{[n]_q !}q^{\frac{n(n-1)}{2}}
  \bigl(e^{\frac{\la}{2} H }X^{+}\oo
  e^{-\frac{\la}{2} H }X^{-}\bigr)^n . \n
\ee
Algebra $U_q(sl(2))$  has three different real forms
determined by their action on the generators as
\begin{displaymath}
\begin{array}{lllrlrll}
\theta_{sl(2,{\bf R})}&:& H \to& H,& X^{\pm} \to& X^{\pm},&
\la \in {i\bf R},\\
\theta_{su(2)} &: & H \to &- H,& X^{\pm} \to & - X^{\mp},&
\la \in  {\bf R}, \\
\theta_{su(1,1)}&:& H \to & - H,& X^{\pm} \to& X^{\mp},&
\la \in  {\bf R},
\end{array}
\end{displaymath}
and homomorphically extended over whole the algebra.

Another possible  solution $U_h(sl(2))$ (we use definition
$h=e^\xi$, which differs from generally accepted)
is called the jordanian quantization and obtained by
twist of the classical algebra  \cite{Ogiv}
with the element
 $\Phi=e^{- \xi H\oo \si}$.
Expression of $\si$ through  $X^-$
is given above where we discuss real forms of the Borel
subalgebra of $sl(2)$. Here we introduce deformation parameter $\xi$
that results in substitution $(X^-, \si) \to (\xi X^-, \xi\si)$
of the previously specified generators.
The universal R-matrix for the jordanian quantization is
\be
{\cal R}_h & = &  \exp\bigl( \xi \si\oo H)
                   \exp\bigl(- \xi H\oo \si).
              \label{RO}
\ee
The only real form
$\tilde\theta_{sl(2,{\bf R})}$  for imaginary $\xi$
is expressed through the classical involution
$H\to H$, $X^{\pm} \to X^{\pm}$ by means of the element
$u=\sum_{n=0}^{\infty}\frac{(-\xi)^n}{n!} H^n \si^n$.
The similarity transformation with the element $e^{\al H}$
turns algebra $U_h(sl(2))$  into $U_{e^{2\al}h}(sl(2))$, so all the
jordanian quantizations with non-zero deformation parameters are isomorphic
to each other.

Classification of the  quantum Lorentz group has been developed
in \cite{PW,WZ}. Below is the list of quantum Lorentz algebras dual
to those solutions of \cite{WZ} which have the classical limit.

\begin{enumerate}
\item
Algebra  $U_{q;r}(so(1,3))$ is built as the twisted tensor product
\twist{$U_q(sl(2))$}{$U_{{\bar q}^{-1}}(sl(2))$} с $\Phi = r^{H_1\oo H_2}$.
For real values of $r$ the real form is defined by the equality
$\theta=\tau(\theta_{sl(2,{\bf R})} \oo \theta_{sl(2,{\bf R})})$,
where $\tau$ denotes the operation of permutation between the
left and right $sl(2)$ subalgebras.
If $q=1$, this solution is  the twist of classical Lorentz algebra
induced by the corresponding deformation of its Cartan subalgebra.
\item
As algebras over the complex field, $U^\pm_{q}(so(1,3))$ are
isomorphic to the twisted tensor square
\twist{$U_{q}(sl(2))$}{$U_{q}(sl(2))$} with
the universal matrix  ${\cal R}_q$  as the twisting element.
For real parameter $q$ there exists two different real forms
$\theta^+=\tau(\theta_{su(2)}\oo \theta_{su(2)})$ and
$\theta^-=\tau(\theta_{su(1,1)} \oo \theta_{su(1,1)})$.
This algebra was introduced in \cite{SWZ,OSWZ}.
\item
Algebra$U_{h;r}(so(1,3))$ is
\twist{$U_{h}(sl(2))$}{$U_{{\bar h}^{-1}}(sl(2))$}
with $\Phi = r^{\si_1 \oo \si_2}$, where dependence on the parameter
$h$ can be eliminated by a change of the basis. As this solutions
is a composition of two twists (the first one is the jordanian
deformation of the left and right $sl(2)$ components),
algebra $U_{h;r}(so(1,3))$ is a pure twist of
the classical Lorentz algebra. Its real form
is given by
$\theta=\tau(\tilde\theta_{sl(2,{\bf R})} \oo \tilde\theta_{sl(2,{\bf R})})$,
for real values of $r$.
\item
Algebra  $U_{h}(so(1,3))$ can be represented as the tensor square of the
jordanian $U_h(sl(2))$ twisted by its R-matrix taken as $\Phi$.
However, the Lorentz real form is easier described if we
consider the direct twist of $U(so(1,3))$ with the element
$\Phi=\exp\{ (H_1+H_2)\oo\frac{1}{2}\ln(1-\xi(X^-_1+X^-_2))\}$.
It is induced by the deformation of the subalgebra isomorphic to
the $sl(2)$ Borel subalgebra, and the real form differs from
classical only by the similarity transformation with the element $u$ as
discussed above. Parameter $\xi$ is assumed to be imaginary.
\end{enumerate}
Thus all the quantum Lorentz algebras belong to the two classes of
twist-equivalence
related to the standard and jordanian deformations. Twist takes part in
quantization of the Poincar\'e algebra as well \cite{LRTN,M7}, although
in this case the explicit description of twist-equivalence classes
is as yet unknown.

\section{Twisted modules}

Twist preserves the multiplicative structure of Hopf algebras
as well as their central elements. Tensor products of representations
of twist-equivalent algebras are isomorphic too, and their spectra
coincide. That does not make any difference between
such Hopf algebras from the point of view of internal
symmetries. Kinematical and dynamical consequences may be rather
significant, however. Quasiclassical limit displays
profoundly unusual behavior of particle trajectories on Minkowski
space even for the simplest deformations \cite{AzKulRod}. So,
throughout this paper one can imagine twist as a transformation
of the classical algebra $U(so(1,3))$.

Recall that a module-algebra ${\cal A}$ for a Hopf algebra
is endowed with a multiplication respected by the action of  ${\cal H}$
in the following sense \cite{ChPr}:
$$
h (a\cdot b)= h^{(1)}a\cdot h^{(2)}b, \quad  h\in {\cal H},
\quad a,b \in {\cal A}.
$$
Here $h^{(1)}\oo h^{(2)}$  denotes symbolically the coproduct  $\Delta(h)$.
This definition matches the algebra of functions on some manifold, or
the differential operator algebra, or just the matrix ring of
linear transformation on a vector space, where  ${\cal H}$
is the enveloping algebra of a classical group of automorphisms.
Twist of  ${\cal H}$ provides a new algebraic structure on
${\cal A}$ consistent with the action of $\tilde {\cal H}$.
This new algebra coinciding with ${\cal A}$ as an ${\cal H}$-module
and denoted further on by $\tilde {\cal A}$, is equipped with
the multiplication
$$
a * b = \Phi_1 a\cdot \Phi_2 b,\quad a,b \in {\cal A}.
$$
Associativity is guaranteed by equation (\ref{TE}). Note that in the case
of function algebra on a classical Poisson manifold this
construction provides an example of deformation quantization introduced in
\cite{BFFLS}.

Properties of two twist-equivalent $q$-manifolds may differ significantly.
Nevertheless, there are transition rules between
geometrical objects if they are somehow related to the symmetry
algebra action. Hardly one may look forward to relate
a representation of $\tilde {\cal A}$ to an arbitrary
representation of ${\cal A}$ as even their spectra may
differ, e.g. due to commutativity violation.
However, this is possible for ${\cal H}$-covariant representations
meaning that both ${\cal A}$ and  ${\cal H}$ act on the same linear space
${\cal W}$, and  $\pi(h a)=\rho(h^{(1)}) \pi(a) \rho(S(h^{(2)}))$
(notations $\pi$ and $\rho$ for  ${\cal A}$- and
${\cal H}$-actions will be omitted further on).
Within the classical
differential geometry such representations correspond to vector
bundles with a group leaf-wise action. In particular, the regular
representation of ${\cal A}$ on itself (trivial 1-bundle) by
multipliers is ${\cal H}$-covariant.  We may also think of
${\cal A}$ as the differential operator algebra on the Minkowski
space and of ${\cal W}$ as the algebra of functions. Let us show
that ${\cal H}$-covariant module ${\cal W}$ over ${\cal A}$ turns
into $\tilde{\cal H}$-covariant module $\widetilde{\cal W}$ over
$\tilde{\cal A}$, if action of $\tilde{\cal A}$ on
$\widetilde{\cal W}$ is  defined as
$$
a*  w = \Phi_1 a \cdot \Phi_2 w,\quad a \in {\cal A},\quad  w \in {\cal W}.
$$
To check this up, it is sufficient to notice that
${\cal H}$-covariance follows (and {\em vice versa}) from
$h (a\cdot  w)= h^{(1)}a\cdot h^{(2)} w$.
The rest implications repeat those for the case ${\cal W} = {\cal A}$.

From geometry of vector bundles, let us proceed to a more detailed
study of the geometry of the base space itself. If an algebra
${\cal A}$ is endowed with an involutive antilinear
anti-automorphism $a \to \bar a$, it can be transferred to
$\tilde{\cal A}$, provided the twisting element satisfies  (\ref{CTE0}).
If that is not the case but identity (\ref{CTE}) holds,
the real form steel can be introduced on
$\tilde{\cal A}$.
First let us prove the following proposition
\newtheorem{lemma}{Lemma}
\begin{lemma}
If the element $\Phi$ satisfies (\ref{CTE}),
the equality
\be
\theta(u)=u^{-1}
\label{uc}
\ee
is true.
\end{lemma}
Using the permutation rule $*S=S^{-1}*$, we find
$$\theta(u)=\theta(\Phi^{-1}_{1})\theta S(\Phi^{-1}_{2}))=
S((\Phi^{-1})^*_{1})S[S(\Phi^{-1}_{2})]^*=
S((\Phi^{-1})^*_{1})(\Phi^{-1})^*_{2}=S(\Phi_{1})\Phi_{2}= u^{-1}.$$
The latter holds for any solution of the twist equation
\cite{D4}. It follows immediately that (\ref{uc})
is true with respect to $\tilde\theta$, too.
We define the twisted real form on the module-algebra
$\tilde{\cal A}$ as $u \bar f$. It is easy to check that
 the operation thus introduced is involutive, anti-homomorphic
 and consistent with the action of $\tilde{\cal H}$.

Suppose now that ${\cal A}$ possesses a measure $\mu$, i.e. a linear
functional positive on elements of the form  $a\cdot \bar a$ (as the
function algebra on a locally compact topological space does).
The same measure is valid for $\tilde{\cal A}$. Indeed,
we find $a* \bar a = \Phi_1 a\cdot \Phi_2 \bar a =
\Phi_1 a\cdot \overline{(\theta({\Phi_2}) a)} $. If
identity (\ref{CTE0}) is fulfilled, the relation
$\Phi_1\oo \theta(\Phi_2) = \theta(\Phi_2) \oo \Phi_1$ holds as well,
and, consequently,  $\Phi_1\oo \theta(\Phi_2)$ can be represented
by a sum $\sum \varphi_i\oo \varphi_i$. Further, we have
$a * \bar a = \sum\varphi_i a\cdot \overline{\varphi_i a}$, and
therefore $\mu(a* \bar a)\geq 0$.
In the case when (\ref{CTE}) is true, one can extend the Hopf algebra
adding the square root of $u$. It is straightforward that the
composition of the coboundary twist with the element
$\Delta(u^{-\frac{1}{2}})(u^\frac{1}{2}\oo u^\frac{1}{2})$ and successive
twist with the element
$(u^{-\frac{1}{2}}\oo u^{-\frac{1}{2}})\Phi (u^\frac{1}{2}\oo u^\frac{1}{2})$
obeys (\ref{CTE0}). This double transformation is carried out
by means of the 2-cocycle
$\Delta(u^{-\frac{1}{2}})\Phi (u^\frac{1}{2}\oo u^\frac{1}{2})$,
and the required property (\ref{CTE0}) readily follows
from (\ref{AuxE}). So, we can apply all the previous considerations
to this composite twist, which differs from initial one by
the internal automorphism of the Hopf algebra only.

By means of $\mu$ an integration and the Hermitean scalar
product are introduced on $\tilde{\cal A}$, turning it
into a pre-Hilbert space. The regular representation of the involutive
algebra $\tilde{\cal A}$ comes out to be a homomorphism into the
$*$-operator algebra  on that space.
It is seen, that although integrals are not changed in  proceeding to
$\tilde{\cal A}$,  the Hermitian form $(a,b)=\mu( a\cdot\bar b)$ is
undergone deformation. Isometries among module-algebras would
remain so among their twisted counterparts, provided they are permutable
with the action of $\tilde{\cal H}$. In particular, this
is true for the classical Fourier transformation, which is a
homomorphism from the algebra of functions on the coordinate
Minkowski space to the convolution algebra on the
space of momenta.

\section{Klein-Gordon-Fock equation}

Taking into account considerations of the previous section,
one can state close connection between structures on
${\cal A}$ and $\tilde{\cal A}$. This conclusion is not of
theoretical interest  only, but yields a machinery
to solve equations of motion, study physical states and so on.
The principal observation is that, providing for some
$c\in {\cal A}$ and all $h\in {\cal H}$ the condition
$h c = \ve(h) c$ is obeyed, the result of the action of $c$
on the elements of the modules ${\cal W}$ and $\tilde{\cal W}$,
coinciding as linear spaces, will be the same:
$c*  w = (\Phi_1 c)\cdot(\Phi_2)  w= (\ve(\Phi_1)c)\cdot (\Phi_2  w)=
c\cdot  w$ for arbitrary $ w\in {\cal W}$.
Analogously, if for some $ w\in {\cal W}$ and any $h\in {\cal H}$
equality $h  w = \ve(h)  w$ takes place, one has $c*  w = c\cdot  w$.
In other words, the result of parings with an invariant
element coincides in the twisted and untwisted cases.
We immediately find that ${\cal H}$-invariant elements of the center of
${\cal A}$ remain so for $\tilde{\cal A}$.
If ${\cal A}$ means the classical Minkowski space ${\cal M}$,
any Lorentz invariant element will commute with the whole
$\tilde{\cal M}$. So, for example, will do the invariant length.

Ordered homogeneous monomials
$(z^{\mu_1}\cdot\ldots \cdot z^{\mu_k})_{k=0,\ldots,\infty}$
in coordinate functions $z^{\mu}$ on the classical Minkowsky
space form a basis of the function algebra ${\cal M}$.
Since the deformed product $*$ is expressed
through $\cdot$ by the invertible element $\Phi$, the quadratic relations
$z^{\mu}z^{\nu}=z^{\nu}z^{\mu}$ in  ${\cal M}$
go over into the quadratic
relations
$(\Phi^{-1}_1 z^{\mu})*(\Phi^{-1}_2 z^{\nu})=(\Phi^{-1}_1 z^{\nu})*(\Phi^{-1}_2 z^{\mu})$,
due to the linear action of ${\cal H}$ on coordinate functions.

Homogeneous monomials of degree $n$ with respect to the twisted product
are expressed through monomials of the same degree in terms of the
old product:
$$ z^{\mu_1}* ... * z^{\mu_k} = \cdot\bigl(\Omega^k(z^{\mu_1}\oo ... \oo  z^{\mu_k})\bigr),$$
where elements $\Omega^k\in {\cal H}^{\oo k}$ form the family of operators
intertwining braid group $B_n$ representations related to the R-matrices
${\cal R}$ and $\tilde{\cal R}$ \cite{KM}. Consequently, ordered monomials
$(z^{\mu_1}*\ldots * z^{\mu_k})_{k=0,\ldots,\infty}$ form a basis
in the twisted module $\tilde{\cal M}$ too, and  functions
on the quantum Minkowski space can be formally expanded over it.
Thus, physically interesting Lorentz-invariant objects remain the same
as elements of twisted modules however they have entirely different expression
through the generators in terms of the twisted multiplication.
Let us turn, for example, to the  Klein-Gordon-Fock equation on the classical
Minkowski space
\be
(g^{\mu\nu}\partial_\mu\cdot\partial_\nu  + m^2) f(z)&=&0.
\label{klein}
\ee
Since the D'Alembert operator is Lorentz-invariant, function $f(z)$,
as belonging to $\tilde{\cal M}$, satisfies the same equation with
D'Alembertian
$$ {\tilde g}^{\mu\nu}\partial_\mu*\partial_\nu=
g^{\mu\nu}(\Phi^{-1})_{\mu\nu}^{\rho\xi}\partial_\rho*\partial_\xi
 = g^{\mu\nu}\partial_\mu\cdot\partial_\nu
.$$
Here  $(\Phi)_{\mu\nu}^{\rho\xi}$  denotes the matrix of the action of
$\Phi$ on partial derivatives $\partial_\mu$. Thus twisting gives
rise to redefinition of the metric. Now consider the function
$e^{i(p,z)}$ on the cotangent bundle over the classical Minkowski
space. The condition
$g^{\mu\nu}p_\mu\cdot p_\nu=m^2$ imposed, it becomes
an element of the algebra of functions  on the mass shell, which is the
quotient of
$\tilde{\cal M}$ over the ideal $(g^{\mu\nu}p_\mu\cdot p_\nu-m^2)$.
Then $e^{i(p,z)}$ is a solution to the quantum Klein-Gordon-Fock
equation and represents a plane wave. Because of the Lorentz-invariance
of $(p,z)$, power series in this element
coincide both in the classical and twisted cases, while the canonical
paring between linear coordinate
space and that of momenta is deformed:
${\tilde g}^\mu_\nu p_\mu* z^\nu={\tilde g}^\mu_\nu
(\Phi_1 p_{\mu})\cdot(\Phi_2 z^{\mu})=g^\mu_\nu p_\mu\cdot z^\nu$.
The general solution to (\ref{klein}) for the real scalar field can be
decomposed into the sum of the plane waves $$ f(z)=\int dp
\delta(p^2-m^2)[a(p)e^{-i(p,z)}+a^\dagger(p)e^{i(p,z)}], $$ that now
can be treated in terms of $\tilde{\cal M}$, with $\int dp
\delta(p^2-m^2)$ being the invariant measure on the quantum mass shell
and the integrand rewritten using the new product $*$.

Let us show how this conclusion agrees with earlier results by
\cite{AzKulRod} taking the simplest case of the algebra
$U_{1;r}(so(1,3))$ as an example. We are interested in the three sets
\begin{displaymath}
Z=\left(\begin{array}{rr}
   z^1  &  z^4   \\
   z^2  &  z^3
\end{array} \right),\quad
P=
\left(\begin{array}{rr}
   p_1  &  p_2   \\
   p_4  &  p_3
\end{array} \right),\quad
D=
\left(\begin{array}{rr}
   \D_1  &  \D_2   \\
   \D_4  &  \D_3
\end{array} \right),
\end{displaymath}
representing, respectively, coordinates, momenta, and  partial derivatives
on ${\cal M}$, which are transformed under the Lorentz group coaction
according to the rule
$$
Z\to M Z M^\dagger,\quad P\to (M^\dagger)^{-1} P M^{-1},
\quad D\to (M^\dagger)^{-1} D M^{-1},\quad M\in SL(2).
$$
Three classical Lorentz invariants are given by the bilinear forms
$$
(z,z)=z^1 z^3 - z^2 z^4,\quad (p,p)=p_1 p_3 - p_2 p_4,\quad
(p,z)=p_1 z^1 + p_2 z^2 + p_3 z^3 + p_4 z^4,
$$
and the matrix commutation relations between $D$ and $Z$
read $D_1 Z_2 - Z_2 D_1 = {\cal P}$.  Here  ${\cal P}$ is the
permutation operator in ${\bf C}^2\oo {\bf C}^2$:
\begin{displaymath}
{\cal P}=
\left(\begin{array}{llll}
  1 & 0 & 0 & 0 \\
  0 & 0 & 1 & 0 \\
  0 & 1 & 0 & 0  \\
  0 & 0 & 0 & 1
\end{array} \right).
\end{displaymath}
Let us perform the twist of the classical universal enveloping
Lorentz algebra with the element $r^{- H_2 \oo H_1}$, where
$H_1$ and $H_2$ are its Cartan generators and consider
the resulting quantum Minkowski space.
Action of  $H_1$ and  $H_2$ on $Z$ and $P$ is defined as
$$
H_1\colon\left(\begin{array}{rr}
   z^1  &  z^4   \\
   z^2  &  z^3
\end{array} \right)\to
\left(\begin{array}{rr}
   z^1   &  z^4   \\
  - z^2  & - z^3
\end{array} \right), \quad
H_1\colon\left(\begin{array}{rr}
   p_1  &  p_2   \\
   p_4  &  p_3
\end{array} \right)\to
\left(\begin{array}{rr}
  - p_1  &   p_2   \\
  - p_4  &   p_3
\end{array} \right),
$$
$$
H_2\colon\left(\begin{array}{rr}
   z^1  &  z^4   \\
   z^2  &  z^3
\end{array} \right)\to
\left(\begin{array}{rr}
   z^1  & - z^4   \\
   z^2  & - z^3
\end{array} \right), \quad
H_2\colon\left(\begin{array}{rr}
   p_1  &  p_2   \\
   p_4  &  p_3
\end{array} \right)\to
\left(\begin{array}{rr}
 - p_1  & - p_2   \\
   p_4  &   p_3
\end{array} \right).
$$
From now on we adopt a convention of distinguishing the
elements of ${\cal M}$, treated as elements of  $\tilde{\cal M}$
by the tilde. Evaluating
multiplication on $\tilde{\cal M}$ we find
$$
\tz^\mu* \tz^\nu = a(\mu,\nu) z^\mu\cdot z^\nu,\quad
\tp_\mu* \tp_\nu = a(\mu,\nu) p_\mu\cdot p_\nu,
$$
$$
\tz^\mu* \tp_\nu = b(\mu,\nu) z^\mu\cdot p_\nu,\quad
\tp_\mu* \tz^\nu = b(\mu,\nu) p_\mu\cdot z^\nu,
$$
where numbers   $a(\mu,\nu)$ and $b(\mu,\nu) = (a(\mu,\nu))^{-1}$
form the matrices (the first index labels the rows)
\begin{displaymath}
a(\mu,\nu) =
\left(\begin{array}{llll}
  r^{-1} & r      & r      & r^{-1} \\
  r^{-1} & r      & r      & r^{-1} \\
  r      & r^{-1} & r^{-1} & r      \\
  r      & r^{-1} & r^{-1} & r
\end{array} \right), \quad
b(\mu,\nu) =
\left(\begin{array}{llll}
  r      & r^{-1} & r^{-1} & r      \\
  r      & r^{-1} & r^{-1} & r      \\
  r^{-1} & r      & r      & r^{-1} \\
  r^{-1} & r      & r      & r^{-1}
\end{array} \right).
\end{displaymath}
This leads to the following relations among
 $\tz$ and $\tp$:
  $$
  \tz^\mu * \tz^\nu = a(\mu,\nu)b(\nu,\mu)\tz^\nu * \tz^\mu
  ,\quad
  \tp_\mu * \tp_\nu = a(\mu,\nu)b(\nu,\mu)\tp_\nu * \tp_\mu,
  $$
  $$
  \tp_\mu * \tz^\nu = b(\mu,\nu)a(\nu,\mu)\tz^\nu * \tp_\mu,
  $$
Explicitly, in the coordinate sector this reads \cite{CD}
\begin{displaymath}
\begin{array}{rrrrrrrrr}
\tz^1 * \tz^2 &=& r^2 \tz^2 * \tz^1, & \tz^1 * \tz^3 &=& \tz^3 * \tz^1, &
\tz^4 * \tz^1 &=& r^2 \tz^1 * \tz^4,\\
\tz^2 * \tz^3 &=& r^2 \tz^3 * \tz^2, & \tz^2 * \tz^4 &=& \tz^4 * \tz^2, &
\tz^3 * \tz^4 &=& r^2 \tz^4 * \tz^3.
\end{array}
\end{displaymath}
Twisting of the Lorentz algebra of the type being considered is
induced by that of its Cartan subalgebra. The latter
preserves two bilinear forms $z_1\cdot z_3$ and $z_2\cdot z_4$, and
that explains why elements  $\tz_1*\tz_3$ and $\tz_2*\tz_4$ belong
to the center of  $\tilde{\cal M}$.

Commutation relations involving partial derivatives $\td_\mu$
are obtained in the similar way:
  $$
  \td_\mu * \tz^\nu = b(\mu,\nu)a(\nu,\mu)\tz^\nu * \td_\mu+b(\mu,\nu)
  ,\quad
  \td_\mu * \tp_\nu = a(\mu,\nu)b(\nu,\mu)\tp_\nu * \td_\mu
  $$
  $$
  \td_\mu * \td_\nu = a(\mu,\nu)b(\nu,\mu)\td_\nu * \td_\mu
  $$
Introducing the operator  $V=r^{\si^3\oo \si^3}$ with
the Pauli matrix  $\si^3=$
$\left(\begin{array}{rr}
   1  &  0    \\
   0  &  -1
\end{array} \right)$,
one can rewrite these relations in the compact form
employed  in \cite{AzKulRod}:
$$\tilde Z_1 V \tilde Z_2 = \tilde Z_2 V \tilde Z_1,\quad
  \tilde D_1 V^{-1} \tilde D_2 = \tilde D_2 V^{-1} \tilde D_1,\quad
  \tilde P_1 V^{-1} \tilde P_2 = \tilde P_2 V^{-1} \tilde P_1,$$
$$\tilde P_1  \tilde Z_2 = V \tilde Z_2  \tilde P_1V^{-1},\quad
  \tilde D_1 V^{-1} \tilde P_2 = \tilde P_2  \tilde D_1.
  $$
The Lorentz invariants  turn into
\be
  (\tz,\tz)_r &=& (z,z)= r^{-1}\tz^1* \tz^3 - r\tz^2* \tz^4, \n\\
  (\tp,\tp)_r &=& (p,p)= r^{-1}\tp_1* \tp_3 - r\tp_2* \tp_4, \n\\
  (\tp,\tz)_r &=& (p,z)= r^{-1}\tp_1* \tz^1 + r\tp_2* \tz^2
                                 + r^{-1}\tp_3* \tz^3 + r\tp_4* \tz^4.\n
\ee
The commutation relations obtained agree with those deduced in
  \cite{AzKulRod}, where the role of $\td_\mu$ is played by
  quantum partial derivatives $\delta_\mu$. However, the matrix
  $\tilde D_1 \tilde Z_2 - V \tilde Z_2\tilde D_1 V^{-1}$
  turns out to be dependent on the deformation parameter
$$
\left(\begin{array}{llll}
  r      & 0      & 0      & 0      \\
  0      & 0      & r^{-1} & 0      \\
  0      & r      & 0      & 0      \\
  0      & 0      & 0      & r^{-1}
\end{array} \right),
$$
in contrast to \cite{AzKulRod}.
The correspondence $\delta_\mu \leftrightarrow \td_\mu$ and interpretation
of the relations involving $\delta_\mu$ in terms of $\partial$ and $\td$
thus require special investigation. Entities $\delta_i$ are introduced
in \cite{AzKulRod} via the transformation law under the
coaction $D\to (M^\dagger)^{-1} D M^{-1}$ of the Lorentz group.
The commutation relations $M^\dagger_1 V M_2=  M_2 V M^\dagger_1$ in the group
determine the invariant uniquely modulo a factor.
In terms of $\delta_\mu$ that invariant is written in the form
$\delta^2= r\delta_1\delta_3 - r^{-1} \delta_2\delta_4$,
whereas the twist of the algebra $U(so(1,3))$ gives
$r^{-1}\tz_1*\tz_3 - r \tz_2*\tz_4$.
It is seen that $\delta^2$ и $(\tp,\tp)_r$ have different dependence
on $r$. On the other hand, deformation of the Lorentz metric
in proceeding from ${\cal M}$ to $\tilde{\cal M}$  implies that the basis
$\td_\mu$ is no more conjugate to the basis
$\tz^\mu$  and matrices  $\td_\mu*\tz^\nu$ и $\tz^\nu*\td_\mu$
behave unsatisfactory under the Lorentz transformations.
For the transition to the  conjugate basis, one should set
$$
(\delta_1,\delta_2,\delta_3,\delta_4)=
(r^{-1}\td_1,r\td_2,r^{-1}\td_3,r\td_4).
$$
Thus one comes to the right expression of the invariant
$r\delta_1\delta_3 - r^{-1} \delta_2\delta_4=r^{-1}\td_1*\td_3
- r \td_2*\td_4.$
The matrix of relations between $\delta_\mu$ and  $\tz^\nu$ takes
the required form independent\footnote{Let us note that there is
a freedom in the choice of the twisting element.
Algebra  $U_{1;r}(so(1,3))$  may be equally obtained if one takes
$\Phi=e^{\al H_1 \oo H_2 + \bt H_2 \oo H_1}$, $\al-\bt=\ln(r)$.
The form of the Lorentz invariant depends a particular choice of
parameters $\al$ and $\bt$. For example,
assuming  $\al=-\bt$, they stay undeformed at all, and two bases
($\tz^\mu$) and ($\tp_\mu$) remain orthonormal and mutually conjugate.}
on $r$.

Deviation from commutativity in $\tilde{\cal M}$ may be accounted for
the so called "undressing" transformation proposed in \cite{AzKulRod}.
It employs two elements $u$ and  $v$, fulfilling the Weyl relations
$vu = r^2 uv$, $ \overline{u} =u^{-1}$, $\overline{v} =v$.
Coordinates and momenta are expressed through $u$ and $v$ and classical
commutative generators
$$
X=
\left(\begin{array}{rr}
    x^1  &   x^4   \\
    x^2  &   x^3
\end{array} \right), \quad
Y=
\left(\begin{array}{rr}
      y_1    &   y_2   \\
      y_4    &   y_3
\end{array} \right)
$$
by the formulas
$$
\tilde Z=
\left(\begin{array}{rr}
   v x^1  &  u^{-1} x^4   \\
   u x^2  &  v^{-1} x^3
\end{array} \right), \quad
\tilde P=
\left(\begin{array}{rr}
   v^{-1} y_1  &  u^{-1} y_4   \\
   u      y_4  &  v      y_3
\end{array} \right). \quad
$$
Noncommutative elements $u$ and $v$ drop from the scalar products
for which we obtain the following:
$$(z,z)=(\tz, \tz)_r = (x,x)_r,\>
  (p,p)=(\tp, \tp)_r = (y,y)_r,\>
  (p,z)=(\tp, \tp)_r = (y,x)_r.$$
These noncommutative elements also disappear from the D'Alembert operator,
hence the function $\int d\tilde\mu(y) [ \tilde a(y)e^{-i(y,x)_r}+h.c.]$
is the  general solution to the quantum Klein-Gordon-Fock equation with
the deformed metric.
Proceeding to the  orthonormal basis shows that this function is equal to
$\tilde f(x)$, где  $f(z)=\int d\mu(p) [ a(p)e^{-i(p,z)}+h.c.]$.
Returning to non-commutative variables in
$\tilde f(x)$,  we find again that $\tilde f(\tilde z)= f(z)$
is the general solution to  the quantum Klein-Gordon-Fock equation.

\section{Dirac equation}

In this paragraph we study four-component spinor fields which
form a ${\cal M}$-module  ${\cal W}={\cal M}\oo {\bf W}$ of
sections of the trivial spinor bundle over the Minkowski space.
The four-dimensional complex lineal ${\bf W}$ is the space of
the $(\frac{1}{2},\frac{1}{2})$-representation of the classical
Lorentz algebra, and the module ${\cal W}$ is freely
 generated by the basic elements $e^{(\al)}\in {\bf W}$:
\begin{displaymath}
e^{(1)}=\left(\begin{array}{c}
             1         \\
             0         \\
             0         \\
             0
\end{array} \right),\quad
e^{(2)}=\left(\begin{array}{c}
             0         \\
             1         \\
             0         \\
             0
\end{array} \right),\quad
e^{(3)}=\left(\begin{array}{c}
             0         \\
             0         \\
             1         \\
             0
\end{array} \right),\quad
e^{(4)}=\left(\begin{array}{c}
             0         \\
             0         \\
             0         \\
             1
\end{array} \right).
\end{displaymath}
The Cartan  generators $H_1$ are  $H_2$ represented by the matrix
\begin{displaymath}
H_1=\left(\begin{array}{rr}
   -\si^3  &  0        \\
     0     &  0
\end{array} \right),\quad
H_2=
\left(\begin{array}{rr}
     0    &    0      \\
     0    &    \si^3
\end{array} \right),
\end{displaymath}
in the chiral basis diagonalizing the matrix $\gm^5$.

The space ${\cal W}$ is also a module over the differential
operator algebra Diff$({\cal W})$, which is the tensor product of
Diff$({\cal M})$  by the ring Lin$({\bf W},{\bf W})$.
In the twisted algebra $\widetilde{\hbox{Diff}}({\cal W})$
linear operators on ${\bf W}$ and partial derivatives do not
commute with each other. Let us calculate the commutation
relations between the Dirac matrices and the generators
$\tz$, $\tp$, and $\td$:
$$
  \tz^\mu * \tgm^\nu = a(\mu,\nu)b(\nu,\mu)\tgm^\nu * \tz^\mu
 ,\quad
  \tp_\mu * \tgm^\nu = b(\mu,\nu)a(\nu,\mu)\tgm^\nu * \tp_\mu,
$$
$$
\td_\mu * \tgm^\nu = b(\mu,\nu)a(\nu,\mu)\tgm^\nu * \td_\mu.
$$
The Dirac matrices themselves obey the relation
$$b(\mu,\nu)\tgm^\mu*\tgm^\nu +b(\nu,\mu)\tgm^\nu*\tgm^\mu =
2 g^{\mu\nu},$$
that gives, in the complex basis of the Minkowski space being
used,
\begin{displaymath}
\begin{array}{rlcrllrll}
  \tgm^1*\tgm^3 +\tgm^3*\tgm^1             &=&  r     ,&
  r^{-1}\tgm^1*\tgm^2  &=& -r     \tgm^2*\tgm^1     ,&
  r     \tgm^1*\tgm^4  &=& -r^{-1}\tgm^4*\tgm^1     ,\\
  \tgm^2*\tgm^4 +\tgm^4*\tgm^2             &=& - r^{-1},&
  r^{-1}\tgm^3*\tgm^4  &=& -r     \tgm^4*\tgm^3     ,&
  r     \tgm^3*\tgm^2  &=& -r^{-1}\tgm^2*\tgm^3     .
\end{array}
\end{displaymath}
Provided permutation rules for the generators and their action on
the basic elements $e^{(\al)}$ are known, it is possible to evaluate
the result of the action of an arbitrary differential operator from
$\widetilde{\hbox{Diff}}({\cal W})$  on an arbitrary element of
$\widetilde{\cal W}$. It is easy to see that
$\td_\mu * e^{(\al)} = 0$.  Further we find
$\tgm^\mu * e^{(\al)} = (\Phi_1\gm^\mu)\cdot (\Phi_2e^{(\al)})$.
In general, the correspondence
$$\tilde{\pi}(\tilde A) = \pi(\Phi_1 A)\rho(\Phi_2) = \rho(\Phi^{(1)}_1) \>\pi(A)\>
\rho(S(\Phi^{(2)}_1)\Phi_2),$$
where we explicitly exposed the representations  $\pi$ and $\rho$ of the
algebras  ${\cal A}$ and ${\cal H}$ on ${\bf W}$, defines a representation
of the twisted module-algebra $\widetilde{\hbox{Lin}}( {\bf W},{\bf W})$
in the ordinary matrix ring  Lin$({\bf W},{\bf W})$,
and that may serve
as another confirmation to the thesis about equivalence between
twisted Hopf algebras as internal symmetries.
Note that the isomorphism mentioned is a homomorphism of the
$\tilde {\cal H}$-modules.
For the Dirac matrices we have
\begin{displaymath}
\tilde{r}(\tgm^1)=
\left(\begin{array}{rrrr}
     0  &  0  &  0 &  0        \\
     0  &  0  &  0 &  1        \\
     r  &  0  &  0 &  0        \\
     0  &  0  &  0 &  0
\end{array} \right),\quad
\tilde{r}(\tgm^2)=
\left(\begin{array}{rrrr}
     0  &  0       &  0 & -1        \\
     0  &  0       &  0 &  0        \\
     0  &  r^{-1}  &  0 &  0        \\
     0  &  0       &  0 &  0
\end{array} \right),\quad
\end{displaymath}
\begin{displaymath}
\tilde{r}(\tgm^3)=
\left(\begin{array}{rrrr}
     0  &  0  &  1 &  0        \\
     0  &  0  &  0 &  0        \\
     0  &  0  &  0 &  0        \\
     0  &  r  &  0 &  0
\end{array} \right),\quad
\tilde{r}(\tgm^4)=
\left(\begin{array}{rrrr}
     0  &  0  &  0 &  0        \\
     0  &  0  & -1 &  0        \\
     0  &  0  &  0 &  0        \\
 r^{-1} &  0  &  0 &  0
\end{array} \right),\quad
\end{displaymath}
These matrices obey the same commutation relations as $\tgm$.

For calculating scattering matrix elements it is necessary to know
how to evaluate traces of twisted matrices. Trace is a
Lorentz-invariant linear functional on ${\hbox{Lin}}( {\bf W},{\bf W})$
which realizes a homomorphism
$\hbox{t}\colon {\cal M} \oo {\hbox{Lin}}( {\bf W},{\bf W}) \to {\cal M}$
between two ${\cal H}$--covariant ${\cal M}$-modules, and
that is true for the twisted case:
$$\tilde{\mbox{t}}(\tilde a*\tilde A) = \mbox{t}(\Phi_1 a\cdot \Phi_2 A)
= (\Phi_1 a) \mbox{t}(\Phi_2 A)=
(\Phi_1 a) \ve(\Phi_2) \mbox{t}(A) =  \tilde a * \tilde{\mbox{t}}(\tilde A)
$$
for $a\in {\cal M}$ and $A\in {\hbox{Lin}}( {\bf W},{\bf W})$.
The functional $\tilde{\mbox{t}}$  can be evaluated
through representation $\tilde \pi$ by the formula
$$\tilde{\mbox{t}} (\tilde A)=\mbox{t} (A) = \hbox{Tr}[r(A)] =
\hbox{Tr}[\tilde \pi (\Phi^{-1}_1 A) \rho (\Phi^{-1}_2)].$$
This expression can be reduced to a more closed expression
held for the $U_{1;r}(so(1,3))$-symmetry:
$\tilde{\mbox{t}} (\tilde A)=\hbox{Tr}[\tilde \pi (\tilde A) \rho (\omega)],$
where $\omega\in {\cal H}$  is equal to
$\Phi^{-1}_1 S(\Phi^{-1}_2) S (\Phi_{2'})\Phi_{1'}$.
Involved in the expression for $\omega$ are the elements taking
part in definition of the twisted antipode $\tilde S$, hence
$\omega=1$ and $\tilde{\mbox{t}} (\tilde A)=\hbox{Tr}[\tilde \pi (\tilde A)]$.
Similar statement is true whenever the twist obeys the factorization
property $(\Delta\oo id)(\Phi)=\Phi_{23}\Phi_{13}$, as that is the case
with the twisted tensor product.

Calculating the trace of the product of two invariant
matrices $(\tp,\tgm)_r$ and $(\tq,\tgm)_r$, where entities $\tq_\mu$
are permuted by the same rules as $\tp_\mu$, gives the invariant scalar
product  $(\tp,\tq)_r$, as might be expected,
having in mind the links between objects of the twist-related
geometries

Now let us turn to the Dirac equation
\begin{equation}
 i \gm^\mu \partial_\mu  \psi(z)  = m \psi(z). \label{Dirac}
\end{equation}
Differential operator  $D=i\gm^\mu \cdot \partial_\mu$  is invariant
under the classical Lorentz algebra action, therefore the quantum
Dirac operator is
$$
   \tilde D=r^{-1}\tgm^1* \td_1 + r\tgm^2 * \td_2
   + r^{-1}\tgm^3*\td_3  + r\tgm^4*\td_4 .
$$
It is preserved by the action of the twisted algebra $U_{1;r}(so(1,3))$.
Eigenvalue functions for the classical Dirac operator in the form of
plane waves read
$$
\begin{array}{llll}
\psi^{(\al)}(z) &=& e^{-i(p,z)}\cdot u^{(\al)}(p),&
\hbox{for positive energies,}\\
\psi^{(\al)}(z) &=& e^{ i(p,z)}\cdot v^{(\al)}(z),&
\hbox{for negative energies,}
\end{array}
$$
where the spinors  $u^{(\al)}(p)$ and $v^{(\al)}(p)$
in the chiral representation are decomposed via the basic elements
$e^{(\al)}$:
\begin{eqnarray}
u^{(1)}(p)&=& (p_1+m)\cdot e^{(1)} + p_4\cdot e^{(2)}
+ (p_3+m)\cdot e^{(3)} -p_4\cdot e^{(4)},\n\\
u^{(2)}(p)&=& p_2\cdot e^{(1)} + (p_3+m)\cdot e^{(2)}
- p_2\cdot e^{(3)} + (p_1+m)\cdot e^{(4)},\n\\
v^{(1)}(p)&=& (p_1+m)\cdot e^{(1)} + p_4\cdot e^{(2)}
- (p_3+m)\cdot e^{(3)} + p_4 \cdot e^{(4)},\n\\
v^{(2)}(p)&=& p_2\cdot e^{(1)} + (p_3+m)\cdot e^{(2)}
+ p_2\cdot e^{(3)} - (p_1+m)\cdot e^{(4)}.\n
\end{eqnarray}
To simplify the resulting expressions we use non-normalized
spinors
$u^{(\al)}(p)$,  $v^{(\al)}(p)$.
Because of the Lorentz-invariance of the Dirac operator, the same
functions expressed through the generators in terms of the new
multiplication will be solutions to the deformed Dirac equation.
The invariant factors $e^{\pm i(p,z)}$  are separated as multipliers
$e^{\pm i(\tp,\tz)_r}$.  Ultimately, we obtain the full set of the
eigenvalue functions of the quantum Dirac operator:
$$
\begin{array}{llll}
\widetilde{\psi}^{(\al)}(\tz) &=& e^{-i(\tp,\tz)_r}*\tu^{(\al)}(\tp),&
 \hbox{for positive energies,}\\
\widetilde{\psi}^{(\al)}(\tz) &=& e^{ i(\tp,\tz)_r}*\tv^{(\al)}(\tp),&
 \hbox{for negative energies}
\end{array}
$$
with
\begin{eqnarray}
\tu^{(1)}(\tp)&=& (r \tp_1+m)*e^{(1)} + r \tp_4*e^{(2)} + (\tp_3+m)*e^{(3)} - \tp_4*e^{(4)},\n\\
\tu^{(2)}(\tp)&=&  r \tp_2*e^{(1)} + (r \tp_3+m)*e^{(2)} - \tp_2*e^{(3)} + (\tp_1+m)*e^{(4)},\n\\
\tv^{(1)}(\tp)&=& (r \tp_1+m)*e^{(1)} + r \tp_4*e^{(2)} - (\tp_3+m)*e^{(3)} + \tp_4*e^{(4)},\n\\
\tv^{(2)}(\tp)&=&  r \tp_2*e^{(1)} + (r \tp_3+m)*e^{(2)} + \tp_2*e^{(3)} - (\tp_1+m)*e^{(4)}. \n
\end{eqnarray}
General solution to (\ref{Dirac})  can be written down as the sum of
the plane waves similarly as for the deformed Klein-Gordon-Fock equation.
It can be derived from the classical solution in two different ways.
So, the transition to the non-commutative variables may be performed
prior to Fourier integration which is then fulfilled via the quantum
invariant measure on the mass shell. Another possible way to obtain a
solution to the deformed Dirac equation from the classical one is to
substitute non-commutative variables into the final classical
expression making use of the connection between the twisted and
non-twisted monomials and bearing in mind that the Lorentz-invariant
factors remain so (as well as the result of multiplication by them)
both in the twisted and  non-twisted cases. This
explains how to deal with the objects of physical interest in proceeding to a
twist-related geometry, no matter classical or already quantum the
original one is.

\section{Conclusion}

In conclusion, let us  briefly formulate the results of our
consideration. We have shown the effectiveness of knowing the
structure of the quantum Lorentz algebra for constructing the
quantum space-time. Group approach appears very natural from
the geometrical point of view, but many nontrivial deformations
are  described in terms of universal enveloping algebras.
In other words, any object of classical geometry transforming
in somehow under the symmetry algebra action does
exist in the quantum space, too, acquiring a new algebraic
content. The problem of description of the algebraic "zoo"
of the phase space thus boils down to
studying unitary representation of the function algebra
on cotangent bundles.

Twists supply with a powerful tool for investigation of the entire
classes of geometries uniformly, and from that point of view it
seems natural to start from the simplest representatives.
In the case of the Lorentz algebra they are $U(so(1,3))$ and
$U_{q;1}(so(1,3))$. Either solutions correspond to
mutually commuting right and left $\frac{1}{2}$- spinors.
Performing twist transformations in accordance with the classification
scheme of Lorentz algebra quantization, one can reach every other
modifications of relativistic geometry, using the advantages
of the presented approach. Rich twist-structure of deformations
of the Lorentz algebra itself makes it possible, in prospective,
to use this technique for building curved spaces of general
relativity, where no translations but only the Lorentz boosts and
rotations survive.

{\bf Acknowledgements.}
One of the authors (PPK) thanks Laurent Baulieu for
the kind invitation in the framework of the jumelage programme,
and the hospitality in LPTHE of Universit\'e Pierre et Marie Curie.

\end{document}